\documentclass{elsarticle}
\usepackage{mathrsfs}
\usepackage{amsmath}
\usepackage{amssymb}
\usepackage{amsfonts}
\usepackage{amsmath,amsfonts,amssymb}

\newtheorem{definition}{Definition}

\newtheorem{theorem}{Theorem}

\newtheorem{proposition}{Proposition}

\newtheorem{lemma}{Lemma}

\newtheorem{remark}{Remark}

\newcommand{\R}{\mathbb R}

\newcommand{\N}{\mathbb N}

\newcommand{\e}{\varepsilon}
\def \O{\Omega}
\def \p{\partial}

\def \l{\lambda}
\def \e{\epsilon}
\def \d{\delta}

\begin{document}
\renewcommand{\thefootnote}{\fnsymbol{footnote}}

\title{Regularity of the extremal solution for some elliptic problems with singular nonlinearity and advection}
\author{Xue Luo$^1$, Dong Ye$^2$, Feng Zhou$^3$}
\address{$^1$Department of Mathematics, East China Normal University,
                       200241 Shanghai, P.R. China. E-mail: luoxue0327@163.com\\
         $^2$LMAM, UMR 7122, Universit\'e Paul Verlaine de Metz, 57045 Metz Cedex 1, France. E-mail: dong.ye@univ-metz.fr\\
         $^3$Department of Mathematics, East China Normal
              University, 200241 Shanghai, P.R. China. E-mail: fzhou@math.ecnu.edu.cn}
\date{}

\begin{abstract}
    In this note, we investigate the regularity of the extremal solution $u^*$
    for the semilinear elliptic equation $-\triangle u+c(x)\cdot\nabla u=\lambda f(u)$
    on a bounded smooth domain of $\mathbb{R}^n$ with
Dirichlet boundary condition. Here $f$ is a positive nondecreasing convex function, exploding at a finite
value $a\in (0, \infty)$. We show that the extremal solution is regular in the low dimensional case. In particular, we prove that for the radial case, all extremal solutions are regular in dimension two.
\end{abstract}


\maketitle

\section{Introduction}
We consider the elliptic problem
\begin{equation}\tag{$P_\lambda$}
    \left\{ \begin{aligned}
         -\triangle u + c(x)\cdot\nabla u & =\lambda f(u) && \textup{in}\
         \Omega,\\
         u& > 0 && \textup{in}\ \Omega,\\
         u&=0 && \textup{on}\ \partial\Omega,
\end{aligned} \right.
\end{equation}
where $\l > 0$, $\Omega$ is a smooth bounded domain in $\mathbb{R}^n$ ($n \geq 2$), $c(x)$
is a smooth vector field over $\overline{\Omega}$ and $f: [0, a) \to \mathbb{R}_+$ with fixed $a \in (0, \infty)$ satisfies the following condition $(H)$:
$$ \mbox{$f$ is $C^2$, positive, nondecreasing and convex in $[0,a)$ with $\displaystyle \lim_{t\to
    a^-}f(t)=\infty$.}$$
In the literature, $f$ is refered as a {\it singular nonlinearity}. We say that $u$ is a regular solution if $u \in C^2(\overline\O)$, and we also deal with solutions in the following weak sense.
\begin{definition}
    We say that $u$ is a weak solution of ($P_\lambda$) if $0\leq
    u\leq a$ a.e. in $\Omega$ such that $f(u)d(x, \p\O)\in L^1(\Omega)$ and
    $$-\int_\Omega u\Delta\phi - \int_\Omega u{\rm div}(\phi c)
    =\lambda\int_\Omega f(u)\phi,
    \quad\forall\; \phi\in C^2(\overline{\Omega})\cap H_0^1(\Omega).$$
Moreover, $u$ is a weak super-solution of ($P_\lambda$) if ${``="}$ is replaced by ${``\geq"}$ for all nonnegative functions $\phi\in C^2(\overline{\Omega})\cap
    H_0^1(\Omega)$.
\end{definition}

Clearly, a weak solution is regular if $\sup_\O u < a$. For regular solutions, we introduce a notion of stability.
\begin{definition}\label{regularity and stability}
    A regular solution $u$ of ($P_\lambda$) is said to be stable if the principal eigenvalue of the
    linearized operator $L_{u,\lambda,c}: =-\triangle + c\cdot\nabla - \lambda f'(u)$
    is nonnegative in $H_0^1(\O)$.
\end{definition}

Exploiting some ideas in \cite{m, gg}, the solvability of $(P_\l)$ is characterized by a parameter $\l^*$:
\begin{proposition}
\label{general}
    There exists $\lambda^*\in(0,\infty)$ such that
\begin{itemize}
 \item For $0<\lambda<\lambda^*$, the problem
($P_\lambda$) has a minimal solution
$u_\lambda$, $u_\l$ is regular and the map $\l \mapsto u_\l$ is increasing. Moreover, $u_\l$ is the unique stable solution of $(P_\l)$.
\item For $\l = \l^*$, $(P_{\l^*})$ admits a unique weak solution $u^* := \lim_{\l\to\l^*} u_\l$, called the extremal solution.
\item For $\l > \l^*$, $(P_\l)$ admits no weak solution.
\end{itemize}
\end{proposition}
Here the minimal solution means that $u_\l \leq v$ for any solution $v$ of $(P_\l)$. We remark immediately a close similarity between $(P_\l)$ and the Emden-Fowler equation with superlinear {\it regular nonlinearity}, that is
\begin{align}
\label{EF}
-\Delta u = \l g(u) \;\mbox{ in } \;\O \subset \R^n; \quad u = 0 \; \mbox{ on }\; \p\O,
\end{align}
with $\l > 0$ and $g: [0, \infty) \rightarrow (0, \infty)$ satisfies
\begin{align}
\label{super}
\mbox{$g$ is $C^2$, nondecreasing, convex and }\; \lim_{t\to \infty} \frac{g(t)}{t} = \infty.
\end{align}
In fact, there exists also a critical parameter $\overline\l \in (0,
\infty)$ for \eqref{EF} such that all conclusions in the above
proposition hold true by replacing $\l^*$ by $\overline\l$ (see \cite{bcmr, m}). It is
well known by classical examples as $g(u) = (1+u)^p$ with $p > 1$ or
$g(u) = e^u$, the extremal solution $u^*$ can be either a regular
solution or a real weak solution in the distribution sense with
$\sup_\O u = \infty$.

\smallskip
For general nonlinearity $g$ satisfying \eqref{super}, the regularity of the extremal solution $u^*$ to \eqref{EF} is obtained by Nedev \cite{n} for any bounded smooth domain $\Omega \subset \mathbb{R}^n$
if $n = 2, 3$; by Cabr\'e \cite{c1} for convex domains in $\R^4$; and for radial symmetry case in $\R^n$ with $n \leq 9$ by Cabr\'e \& Capella \cite{cc}. In \cite{yz1}, it is proved that, under mild condition on $g$, the extremal solution $u^*$ is regular for any smooth bounded domain $\O \subset \R^n$
if $n \leq 9$.

\medskip
We can ask the same question about the problem ($P_\lambda$): For $f$ verifying $(H)$, is it true that the extremal solution to ($P_\lambda$) is regular for general vector field $c$ and general domain $\Omega \subset\mathbb{R}^n$ with low dimensions $n$? We will partly answer this question. It is worthy to mention that for studying the explosion phenomena in a flow, Berestycki {\it et al.} \cite{bknr} have considered the problem $(P_\l)$ with a general source $f$ verifying \eqref{super}.

\medskip
Without loss of generality, fix $a = 1$ in the sequel. The problem $(P_\l)$ can be linked to equation \eqref{EF} up to the transformation $v = -\ln{(1-u)}$. In fact, let $u$ solve $(P_\l)$, $v$ verifies then
\begin{equation}\tag{$Q_\lambda$}
   \left\{ \begin{aligned}
         -\triangle v+|\nabla v|^2+c(x)\cdot\nabla v&=\lambda e^vf(1-e^{-v}):=\lambda g(v) && \textup{in}\ \Omega,\\
         v&=0 && \textup{on}\ \partial\Omega.
\end{aligned} \right.
\end{equation}
Therefore $g$ verifies \eqref{super} and $v^*=-\ln(1-u^*)$ is the extremal solution for the problem ($Q_\lambda$). Thus the regularity of $u^*$ is equivalent to the boundedness of $v^*$, however the situation could be very different with the presence of advection terms (see \cite{cg, wy}). In last decade, a model describing the steady state of MEMS (Micro-Electro-Mechanical Systems) device given by Pelesko and Bernstein in \cite{pb}, has drawn many attentions (see \cite{egg} and the references therein).
$$-\Delta u = \frac{\l}{(1 - u)^2} \;\mbox{ in } \;\O \subset \R^n; \quad u = 0 \; \mbox{ on }\; \p\O.$$
More generally, many precise studies have been done for the singular nonlinearities with negative exponent $f(u) = (1 -u)^{-p}$ ($p > 0$) in the advection-free situation, i.e. $c \equiv 0$.
In that case, when $\O $ is moreover the unit ball in $\R^n$, it is known that $u^*$ is regular if and only if (see \cite{mp, gg})
\begin{align}
\label{np}
n < n_p:= 2 + \frac{4p}{p+1} + 4\sqrt{\frac{p}{p+1}}.\end{align}
Tending $p \to 0^+$ in \eqref{np}, we see that $n_p \to 2$. Therefore we cannot expect in general better than dimension two to claim the regularity of $u^*$.

\medskip
For the radial case of $(P_\l)$, equally when $\O$ is a ball and $c(x)$ is the gradient of a smooth radial function, $u_\l$ is radial by uniqueness of the minimal solution. We obtain the following optimal results which are new even for the advection-free case.
\begin{theorem}
\label{reguradial}
Assume that $n = 2$, $\O = B_1$. Let $\gamma$ is a smooth radial function and $c = \nabla\gamma$, then the extremal solution $u^*$ is regular for any $f$ satisfying $(H)$.
\end{theorem}

\begin{theorem}
\label{radialu'}
For any $f$ satisfying $(H)$, $\O = B_1$ and smooth radial function $\gamma$, there exists $C > 0$ such that for all $\l \in (0, \l^*]$
$$|u_\l'(r)| \leq \left\{\begin{array}{ll}
Cr^{-1} & \mbox{if } n \geq 10;\\
Cr^{-\frac{n}{2}+1+\sqrt{n-1}} & \mbox{if } 3 \leq n \leq 9;
\end{array}
\right.  \quad \forall \; r = {|x|} \in (0, 1]$$
where ${|\cdot|}$ is the Euclidean norm in $\R^n$.
\end{theorem}
\begin{remark}
\label{1.1}
The above estimates are optimal. In fact, when $f(u) = (1 - u)^{-p}$, $p > 0$, $\O = B_1$ and $c \equiv 0$, it is well known that $u^*(x) = 1 - r^{\frac{2}{p+1}}$ if $n \geq n_p$ with $n_p$ given in \eqref{np}, and we have
$$n \geq n_p \quad \mbox{ iff }\quad n \geq 10 \;\mbox{ or }\; 3\leq n \leq 9, \; \frac{2}{p+1} \leq -\frac{n}{2} + 2 + \sqrt{n - 1}.$$
\end{remark}

But is the extremal solution $u^*$ of $(P_\l)$ regular with general singular nonlinearity $f$ verifying $(H)$, vector field $c$ and smooth bounded domains in $\R^2$? The answer is affirmative under some additional mild condition on $f$.
\begin{theorem}\label{main theorem}
Assume that $f$ satisfies conditions $(H)$ and the additional conditions,
$$
            \limsup_{t\rightarrow1^-}\frac {f(t)}{f'(t)(1-t)\ln^2(1-t)}<1 \leqno{(H1)}$$
and
$$\liminf_{t\rightarrow1^-} \frac
            {f(t)f''(t)}{f'^2(t)} > 0.\leqno{(H2)}$$
Then $u^*$ is regular solution to ($P_{\l^*}$) if $n = 2$, i.e.
    $\Omega\subset\mathbb{R}^2$.
\end{theorem}

Under more precise conditions on the growth of $f$, the
extremal solution can be showed to be regular in some higher dimensions.
\begin{theorem}\label{main result}
    Let $f$ verify (H) and $g(v) = e^vf(1 - e^{-v})$. Assume that $g$ satisfies
$$\liminf_{t\rightarrow\infty}
        \frac{g'(t)}{g(t)} = 1+\delta > 1\leqno{(H3)}$$
and
$$\liminf_{t\rightarrow\infty} \frac {g''(t)g(t)}
        {g'^2(t)}=\mu>\frac 1 {1+\delta}. \leqno{(\widetilde{H2})}$$
Then $v^* = -\ln(1 - u^*)$ is bounded (so $u^*$ is regular) when
\begin{align}
\label{ndelta}
n < 2+\frac {4\delta}{1+\delta}+\frac
    {4\sqrt{\delta(\mu + \mu\delta -1)}}{1+\delta}.
    \end{align}
\end{theorem}

Consequently, if $\mu\delta > 1$, $u^*$ is regular for all $n\leq 6$. Furthermore, if we can tend $\delta$ to $\infty$, which means $g = o(g')$ near $\infty$, then $u^*$ is regular for $n < 6 + 4\sqrt{\mu}$ with any $\mu > 0$. However, we can never have $\mu > 1$, since otherwise $g$ blows up at finite value and contradicts \eqref{super}, so the best result we can expect is for $n \leq 9$. For example, if $f(u) = e^{\frac{1}{1 - u}}$, then $g(v) = e^{v + e^v}$ verifies $\delta = \infty$ and $\mu = 1$.

\begin{theorem}\label{fast increasing nonlinearity}
    Let $f$ verify (H) and $g(v) = e^vf(1 - e^{-v})$. Assume that $g = o(g')$ near $\infty$. Rewrite $g(t) = g(0) + te^{h(t)}$ in $(0, \infty)$, suppose there exists $t_0>0$ such that $t^2h'(t)$ is nondecreasing
    for $t \geq t_0$, then for any bounded smooth domain $\Omega\subset \mathbb{R}^n$
    with $n\leq 9$, $u^*$ is a regular solution.
\end{theorem}

Furthermore, when $g = o(g')$ near $\infty$, the condition $(\widetilde{H2})$ is just equivalent to $(H2)$, since
$$\frac{f''(t)f(t)}{f'^2(t)} = \frac{\left(g'' - g'\right)g}{\left(g' - g\right)^2}(s) = \left(\frac{g''g}{g'^2} - \frac{g}{g'}\right)\times \left(1 - \frac{g}{g'}\right)^{-2}(s), \qquad \forall\; t = 1 - e^{-s}.$$
It is also easy to see that $(H3)$ is equivalent to the condition
$$\liminf_{t\rightarrow 1^-}
        \frac{f'(t)(1-t)}{f(t)} = \delta > 0.$$
If the equality holds for the whole limit, we have the following optimal result. The case $f(u) = (1 - u)^{-2}$ was obtained in \cite{cg} with a different argument.
\begin{theorem}\label{negative p}
Assume that
\begin{align}
\label{H2'}
\lim_{u\to 1^-}\frac{f'(u)(1 - u)}{f(u)} = p > 0.
\end{align}
Then $u^*$ is a regular solution if $n < n_p$ where $n_p$ is defined in \eqref{np}.
\end{theorem}

One of the main difficulties here is due to the vector field $c(x)$. When $c \ne 0$, the operator $-\Delta + c\cdot\nabla$ is not self-adjoint, we use ideas from \cite{cg} to get some energy estimates. However if $c$ is a gradient, say $c=-\nabla\gamma$ in $\O$, then $-\Delta +c\cdot\nabla$ can be rewritten as
$e^{-\gamma}L_\gamma$ where $L_\gamma = -{\rm div}(e^\gamma\nabla)$ is a self-adjoint operator. In that case, $(P_\l)$ admits a variational structure and we can expect more precise estimates of minimal solutions $u_\l$, as in the radial case.

\medskip
The paper is organized as follows: In section 2, we prove quickly Proposition \ref{general} and show some general consequences of the stability of $u_\l$. The section 3 is devoted to the proof of Theorems \ref{main theorem} to \ref{negative p} for general domains. In section 4, we discuss the radial case. The norm $\|\cdot\|_q$ denotes always the standard $L^q$ norm for any $q \in [1, \infty]$. The capital letter $C$ denotes a generic positive constant independent of $\l$, it could be changed from one line to another.

\section{Preliminaries}
\setcounter{equation}{0}
\setcounter{remark}{0}
As mentioned above, $-\Delta + c\cdot\nabla$ is not a self-adjoint operator for general vector field $c$. However using Lemma 1 in \cite{cg}, we have a kind of Hodge decomposition, which tells us that for any vector field $c \in C^\infty(\overline{\Omega},\mathbb{R}^n)$, there exist a smooth scalar function $\gamma$ and a vector field $b \in C^\infty(\overline{\Omega},\mathbb{R}^n)$ such that
\begin{align}
\label{decomp-c}
c = -\nabla\gamma+b \; \mbox{ and }\; {\rm div}(e^\gamma b)= 0 \quad \mbox{in }\; \overline\O.
\end{align}
Therefore the problem ($P_\lambda$) can be
rewritten as
\begin{equation}\tag{$P_\lambda'$}
    -{\rm div}(e^\gamma\nabla u)+ e^\gamma b\cdot\nabla u=\lambda
    e^\gamma f(u)\quad\textup{in}\ \Omega.
\end{equation}
On the other hand, we don't have a
suitable variational characterization in general to use the stability
assumption. Fortunately, we can adopt an energy inequality as in \cite{cg}, which is
derived from a generalized Hardy inequality of \cite{co}.
\begin{proposition}
    Let $u_\l$ be minimal solution of $(P_\l)$. For any $1\leq\beta<2$, we have
    \begin{align}\label{Hardy ineqn}
        \lambda\int_\Omega e^\gamma f'(u_\lambda)\psi^2
            \leq \frac2\beta\int_\Omega
            e^\gamma|\nabla\psi|^2+\frac{\|b\|_\infty^2}{2(2-\beta)}\int_\Omega
            e^\gamma\psi^2, \quad \forall\; \psi\in H_0^1(\Omega).
    \end{align}
    where $b$ is the vector field in \eqref{decomp-c}, $\|b\|_\infty = \max_{\overline\Omega}|b(x)|$.
\end{proposition}
\noindent{\bf Proof.} We use a Hardy type inequality given by Theorem 2 in \cite{cg}, which says
that for a positive principal eigenfunction $\varphi$ of $L_{u_\lambda, \l, c}$, for $\beta \in [1, 2)$ and any $\psi \in H_0^1(\O)$,
    $$\lambda\int_\Omega e^\gamma f'(u_\lambda)\psi^2\leq\frac2\beta\int_\Omega e^\gamma|\nabla\psi|^2
        +\int_\Omega \left[-\frac{2-\beta}2\frac{|\nabla\varphi|^2}{\varphi^2}+
        \frac{b\cdot\nabla\varphi}\varphi\right]e^\gamma\psi^2.$$
    By Cauchy-Schwarz inequality, it is easy to see $$-\frac{2-\beta}2\frac{|\nabla\varphi|^2}{\varphi^2}+\frac{b\cdot\nabla\varphi}\varphi \leq \frac{|b(x)|^2}{2(2-\beta)} \leq \frac{\|b\|_\infty^2}{2(2-\beta)},$$
    so we are done.\hfill{$\Box$}

\medskip
Another main ingredient of our approach is just the transformation $v = -\ln(1 - u)$. Let $\phi$ and $\xi$ be nonnegative $C^1$ functions satisfying
$\phi(0) = \xi(0) = 0$ and $\xi'=\phi'^2$. Define $v_\l = -\ln(1 - u_\l)$ and $g(v_\l) = e^{v_\l}f(1 - e^{-{v_\l}})$. Using $(Q_\l)$, we get $-{\rm div}(e^\gamma\nabla v_\l) + e^\gamma b\cdot\nabla v_\l \leq \l e^\gamma g(v_\l)$ in $\O$. Let $\psi = \phi(v_\l)$ in (\ref{Hardy ineqn}), $\forall\; \l \in (0, \l^*)$,
\begin{align*}
    &\; \lambda\int_\Omega e^\gamma
    f'(u_\lambda)\phi^2(v_\lambda)\\
    \leq &\;\frac2\beta\int_\Omega e^\gamma|\nabla\phi(v_\lambda)|^2
    +\frac{\|b\|_\infty^2}{2(2-\beta)}\int_\Omega e^\gamma\phi^2(v_\lambda)\\
    = &\; \frac2\beta\int_\Omega e^\gamma\nabla \xi(v_\lambda) \nabla
    v_\lambda + C_\beta\int_\Omega e^\gamma\phi^2(v_\lambda)\\
    = &\; -\frac2\beta\int_\Omega {\rm div}(e^\gamma\nabla
    v_\lambda)\xi(v_\lambda)+ C_\beta\int_\Omega e^\gamma\phi^2(v_\lambda)\\
    \leq &\; \frac{2\lambda}\beta\int_\Omega e^\gamma g(v_\lambda)\xi(v_\lambda) - \frac{2}\beta\int_\Omega e^\gamma b\cdot \xi(v_\l)\nabla v_\l + C_\beta\int_\Omega
    e^\gamma\phi^2(v_\lambda)\\
= &\; \frac{2\lambda}\beta\int_\Omega e^\gamma g(v_\lambda)\xi(v_\lambda) + C_\beta\int_\Omega
    e^\gamma\phi^2(v_\lambda).
\end{align*}
The last line is due to ${\rm div}(e^\gamma b) = 0$. We claim then
\begin{proposition}
Let $1\leq \beta < 2$. For any $\l \in (0, \l^*)$ and any nonnegative $C^1$ test functions $\phi$, $\xi$ verifying
$\phi(0) = \xi(0) = 0$ and $\xi'=\phi'^2$, there hold
\begin{align}
\label{psi and xi}
\lambda\int_\Omega e^\gamma f'(u_\lambda)\phi^2(v_\lambda)
    \leq \frac{2\lambda}\beta\int_\Omega e^\gamma g(v_\lambda)\xi(v_\lambda) + C_\beta\int_\Omega
    e^\gamma\phi^2(v_\lambda)
\end{align}
and
\begin{align}
\label{est1} \lambda\int_\Omega e^\gamma f'(u_\l)\phi^2(u_\l) \leq
\frac{2\lambda}{\beta} \int_\Omega e^\gamma f(u_\l)\xi(u_\l) + C_\beta
\int_\Omega e^\gamma \phi^2(u_\l).
\end{align}
\end{proposition}
The proof of \eqref{est1} is completely similar to \eqref{psi and xi} but using $(P'_\l)$ instead of $(Q_\l)$.

\medskip
We also make use the following behavior of $f$ proved in \cite{yz2}.
\begin{lemma}
\label{ff'} For any $f$ verifying $(H)$, we have $\lim_{t\to 1} f(t)/f'(t) = 0$.
\end{lemma}

Choose first $\phi(u) = e^u - 1$ in \eqref{est1}, then $\xi(u) = \frac{e^{2u} - 1}{2}$ and
\begin{align*}
\lambda\int_\Omega e^\gamma f'(u_\l)\left(e^{u_\l} - 1\right)^2 \leq
\frac{\lambda}{\beta} \int_\Omega e^\gamma f(u_\l)\left(e^{2u_\l} - 1\right) + C_\beta
\int_\Omega e^\gamma \left(e^{u_\l} - 1\right)^2.
\end{align*}
Fix $\beta \in (1, 2)$. By Lemma \ref{ff'},
$$\lambda\int_\Omega e^\gamma f'(u_\l)e^{2u_\l} \leq C.$$
Consequently $\|f'(u_\l)\|_1$ is uniformly bounded, so is $\|f(u_\l)\|_1$. Multiplying $(P_\l)$ by $u_\l$,
$$\int_\O |\nabla u_\l|^2 = \int_\O \frac{{\rm div}(c)}{2}u_\l^2 + \l\int_\O f(u_\l)u_\l \leq C,$$
which gives
\begin{proposition}
\label{H01}
The family of minimal solutions $\{u_\l\}_{0 < \l < \l^*}$ is uniformly bounded in $H^1_0(\O)$.
\end{proposition}
\begin{remark}
As far as we know, it is always an open question whether the similar $H^1$ energy estimation holds for minimal solutions of \eqref{EF} with general regular nonlinearity satisfying \eqref{super} and general domain $\O$ when $n \geq 6$ (see \cite{n} for $n\leq 5$). For the advection-free case $c = 0$, it was proved in \cite{yz2} that $u^* \in H^2\cap H_0^1(\O)$ under the condition $(H)$, it is also true for the gradient case $c = \nabla \gamma$ (see Lemma \ref{gradientH2}).
\end{remark}

\noindent
{\bf Sketches of proof of Proposition \ref{general}.} We follow the ideas coming from \cite{bknr, m, gg}. The main argument is the maximum principle for operators $-\Delta + c\cdot\nabla$ and $L_\gamma$ under the Dirichlet boundary condition, we use also the super-sub solution method and monotone iteration.

\medskip
Let $w \in H_0^1(\O)$ be the regular solution of $-\Delta w + c\cdot\nabla w = 1$ in $\O$ and fix $\alpha > 0$ such that $\alpha\max_\O w < 1$. It is easy to verify that $\alpha w$ is a supersolution of $(P_\l)$ for $\l > 0$ small enough. As $0$ is a subsolution and $\alpha w > 0$ in $\O$, $(P_\l)$ admits a regular solution for $\l > 0$ small enough. As any regular solution $u$ of $(P_\l)$ is also a supersolution for $(P_\mu)$ if $\mu \in (0, \l)$, the set of $\l$ for which $(P_\l)$ admits a regular solution is just an interval. Moreover, for these $\l$, using $(H)$ and the monotone iteration $v_0 = 0$; $-\Delta v_{n+1} + c\cdot\nabla v_{n+1} = \l f(v_n)$ in $\O$ with $v_{n+1} = 0$ on $\p\O$ for $n \in \N$, we get the minimal solution $u_\l = \lim_{n\to\infty} v_n$.

\smallskip
If we suppose that the principal eigenvalue of $L_{u_\l, \l, c}$ is negative, we can construct, as in \cite{bknr} another solution $v \leq u_\l$ using the associated first eigenfunction, this is just impossible by the definition of $u_\l$, hence $u_\l$ is stable. The uniqueness of stable solution comes from Lemmas 2.16 and 2.17 in \cite{CR}.

\medskip
Take a positive first eigenfunction $\varphi$ of $L_\gamma$ with the Dirichlet boundary condition, by $(P_\l')$,
$$\l f(0)\int_\O e^\gamma \varphi \leq \int_\O\l e^\gamma f(u)\varphi = \int_\O {\l_1(L_\gamma)}u\varphi - \int_\O {\rm div}(e^\gamma b \varphi)u \leq C.$$
So $\l$ is upper bounded. Define the critical threshold $\lambda^*$ as the supermum of $\lambda > 0$
for which ($P_\lambda$) admits a regular solution, as $u^*$ is the monotone limit of $u_\l$ when $\l \to \l^*$,  we deduce that $u^* \in H^1_0(\O)$ is a weak solution of $(P_\l)$ by Proposition \ref{H01}.

\medskip
Suppose that $u$ is a weak solution to $(P_\l)$. By the monotonicity of $f$, it is easy to verify that for any $\delta > 1$, the function $v = \delta^{-1}u$ is a weak supersolution for $(P_{\l/\delta})$, then the monotone iteration will enable us a weak solution $w$ of $(P_{\l/\delta})$ satisfying
$0 \leq w \leq v \leq \delta^{-1} < 1$. The regularity theory implies then $w$ is a regular solution of $(P_{\l/\delta})$. This means that $\l/\delta \leq \l^*$. Let $\delta$ tend to 1, we get $\l \leq \l^*$. Therefore, no weak solution exists for $\l > \l^*$.

\medskip
The uniqueness of the weak solution can be proved in the very similar way as in \cite{m} using the monotonicity and convexity of $f$, with the strong maximum principle for the operator $-\Delta + c\cdot\nabla$ associated to Dirichlet boundary condition, so we omit the details.\hfill{$\Box$}

\section{Regularity of $u^*$ for general $c$ and $\O$}
\setcounter{equation}{0}
\setcounter{lemma}{0}

For proving our results, we will choose suitable functions $\phi$ to apply \eqref{psi and xi} or \eqref{est1}. We need also
\begin{lemma}
\label{comp}
For any $q > n/2$, there exists $C > 0$ such that the solution $v$ of $(Q_\l)$ satisfies $0 \leq v \leq C\|g(v)\|_q$ in $\O$.
\end{lemma}
Indeed, let $w$ be the solution of $L(w) := -\Delta w + c\cdot\nabla w = \lambda g(v)$ in $\O$
with $w = 0$ on $\p\O$. By regularity theory and Sobolev embedding, $\|w\|_\infty \leq C\|w\|_{W^{2, q}(\O)} \leq C'\lambda^*\|g(v)\|_q$ because $q > n/2 \geq 1$. Morover, as $L(w - v) \geq 0$, the maximum principle implies then $0 \leq v \leq w \leq C\|g(v)\|_q$.

\subsection{Proof of Theorem 1.3}
For simplicity, we omit the index $\l$ for $u_\l$ or $v_\l$. Let $\phi(u) = v = -\ln(1 - u)$ in \eqref{est1}, so $\xi(u) = (1 - u)^{-1} - 1$.
Fix $\beta \in (1, 2)$ but very close to 2. Repeating the proof of Theorem 2 in
\cite{yz2} with the assumption $(H1)$, there exists $C > 0$ such that
$$\lambda\int_\Omega e^\gamma \frac{f(u)}{1 - u} < C + C C_\beta\int_\Omega e^\gamma \phi^2(u).$$
As $\phi^2(u) = o(\xi(u)) = o(f\xi)$ when $u \to 1^-$,
$$\lambda\int_\Omega e^\gamma \frac{f(u)}{1 - u} \leq C.$$
Using the equation $(Q_\l)$ and $\p_\nu v \leq 0$ on $\p\O$,
\begin{align*}
\int_\Omega |\nabla v|^2 = \lambda\int e^v f(1 - e^{-v}) +
\int_{\partial\Omega} \frac{\partial v}{\partial \nu} d\sigma
-\int_\Omega c\cdot\nabla v &\leq \lambda\int_\Omega \frac{f(u)}{1 - u} + C\|\nabla v\|_2\\
 & \leq C + C\|\nabla v\|_2.
\end{align*}
Therefore $\|\nabla v\|_2 \leq C$, the classical Moser-Trudinger inequality enables us, as $n = 2$
\begin{align}
\label{mt}
\int_\Omega e^{qv} \leq C_q, \quad \forall\; q \geq 1.
\end{align}
Take now $\phi(u) = f(u) - f(0)$ in \eqref{est1}, we
need to estimate
\begin{align*}
\zeta(u) := f'(u)\phi(u) - \frac{2}{\beta}\xi(u) &= f'(u)\phi(u) - \frac{2}{\beta}\int_0^u f'^2(s) ds\\
&= f'(u)f(u) -
\frac{2}{\beta}\int_0^u f'^2(s) ds - Cf'(u)\\
 &:= I(u) -
\frac{2}{\beta}J(u) - Cf'(u).
\end{align*}
 By $(H2)$, there exists $\delta
> 0$ such that
\begin{align*}
I(u) - I(0) = \int_0^u \left[f'^2(s) + f''(s)f(s)\right] ds \geq (1 +
\delta)J(u) - Cf'(u), \;\; \forall\; u \in [0, 1)
\end{align*}
Let $\frac{4}{2 + \delta} < \beta < 2$, we get $\zeta(u) \geq C I(u) - C$. Asserting this in \eqref{est1},
$$\lambda\int_\Omega e^\gamma f'(u)f^2(u) \leq C\int_\Omega e^\gamma f^2(u) + C.$$
Consequently, $\|f'(u)f^2(u)\|_1 \leq C$. By Lemma \ref{ff'}, we deduce $\|f(u)\|_3 \leq C$. Combining with \eqref{mt}, $\|g(v)\|_p \leq C$ for any $p < 3$. The proof is completed by Lemma \ref{comp} as $n = 2$. \hfill{$\Box$}

\subsection{Proof of Theorem 1.4}
Without loss of generality, we can assume that $g(0)=1$. Let $\phi(t)= g^\alpha(t)-1$
where $\alpha > 0$ is a constant to be determined later. Then
\begin{align}
\label{3.5}
\begin{split}
    \xi(t) & =\int_0^t \phi'^2(s)ds\\
    & =\alpha^2\int_0^t g^{2\alpha-2}(s)g'^2(s)ds\\
    &=\frac {\alpha^2}{2\alpha-1}g^{2\alpha-1}(t)g'(t)-\frac {\alpha^2}{2\alpha-1}\int_0^t
    g^{2\alpha-1}(s)g''(s)ds - C_\alpha.
    \end{split}
\end{align}
The condition $(\widetilde{H2})$ yields: Given any $\epsilon \in \left(0, \mu-\frac 1{1+\delta}\right)$, there exists
$C\geq 0$ such that $g(t)g''(t)\geq(\mu-\epsilon)g'^2(t)-C$ in $[0, \infty)$. Therefore
\begin{align}
\label{3.6}
\begin{split}
    -\int_0^t g^{2\alpha-1}(s)g''(s)ds
    &\leq -(\mu-\epsilon)\int_0^t
    g^{2\alpha-2}(s)g'^2(s)ds+C\\
& \leq -\frac {\mu-\epsilon}{\alpha^2}\xi(t)+C.
\end{split}
\end{align}
We divide the proof into two cases.

\medskip\noindent
{\sl Case 1}: $\delta > 1$ and $\mu > \frac{1}{1 + \delta}$; or $\delta \leq 1$ with $\mu > \frac{1+\delta}{4\delta}$.

\smallskip
Take $\alpha > \frac{1}{2}$. Combine \eqref{3.5} and \eqref{3.6},
$$\left(1+\frac {\mu-\epsilon}{2\alpha-1}\right)\xi(t)\leq\frac
{\alpha^2}{2\alpha-1}g^{2\alpha-1}(t)g'(t)+C,$$
consequently
\begin{align}
\label{xi-g}
\xi(t)\leq \frac
{\alpha^2}{2\alpha-1+\mu-\epsilon}g^{2\alpha-1}(t)g'(t)+C, \quad \mbox{for any }\; t\geq0.
\end{align}
According to $(H3)$, for any $0 < \delta' < \delta$, there exists $C > 0$ such that $g'(t) \geq (1 + \delta')g(t) - C$ in $[0, \infty)$. Setting these estimates in (\ref{psi and xi}), omitting the index $\l$ and recalling that $f'(u) = g'(v) - g(v)$,
\begin{align*}
 & \frac{\d'\l}{1+\delta'}\int_\Omega e^\gamma g'(v)(g^\alpha(v)-1)^2 - C\l\int_\Omega
    e^\gamma(g^\alpha(v)-1)^2\\
    \leq & \; \l\int_\Omega e^\gamma f'(u)(g^\alpha(v)-1)^2\\
    \leq & \;\frac {2\alpha^2\l}{\beta(2\alpha-1+\mu-\epsilon)}\int_\Omega
    e^\gamma g^{2\alpha}(v)g'(v)+ C\l\int_\Omega
    e^\gamma g(v) + C\int_\Omega e^\gamma (g^\alpha(v)-1)^2.
\end{align*}
Consequently,
\begin{align*}
& \left[\frac{\d'}{1+\d'}-\frac
    {2\alpha^2}{\beta(2\alpha-1+\mu-\epsilon)}\right]\l\int_\Omega
    e^\gamma g'(v)g^{2\alpha}(v)\\
\leq & \;\frac {2\d'C}{1+\d'}\int_\Omega
    e^\gamma g'(v)g^\alpha(v) +C\int_\Omega
    e^\gamma g(v) + C\int_\Omega e^\gamma (g^\alpha(v)-1)^2.
\end{align*}
Choose $\d'$ near $\d$ such that
$$\mbox{either } \; \d' > 1 \;\mbox{ and }\; \mu > \frac{1}{1 + \d'} \quad \mbox{or} \quad \d' < \delta \leq 1 \;\mbox{ with } \;\mu > \frac{1+\d'}{4\d'}.$$
Through direct computations, for $\e > 0$ sufficiently small and $\beta = 2 - \e$, there exists $$\alpha \in \left(\frac12,
\frac{\d'}{1+\d'} + \frac{\sqrt{\d'(1+\d')(\mu-\epsilon)-\d'}}{1+\d'}\right)$$
such that
\begin{align}
\label{da}\left[\frac{\d'}{1+\d'}-\frac
    {2\alpha^2}{\beta(2\alpha-1+\mu-\epsilon)}\right] > 0.\end{align}
For such $\alpha$, we obtain
\begin{equation}\label{uniform control}
    \l\int_\Omega e^\gamma g^{2\alpha}(v)g'(v)\leq C, \quad \forall\; \l \in (0, \l^*).
\end{equation}
Tending now $\d'$ to $\d$ and $\e$ to $0$, \eqref{uniform control} holds true provided that
\begin{align}
\label{esta}
\alpha < \frac
\delta{1+\delta}+\frac
{\sqrt{\delta\mu(1+\delta)-\d}}{1+\delta}.\end{align}
Therefore
\begin{align*}
    \int_\Omega e^\gamma g^{2\alpha+1}(v)\leq C\int_\Omega
    e^\gamma g^{2\alpha}(v)g'(v) + C \leq \widetilde C,
\end{align*}
which implies that $\|g(v)\|_{2\alpha+1} \leq C$ for $\alpha$ verifying \eqref{esta}.
Applying Lemma \ref{comp}, we conclude that
for $n < 2 + 4\alpha$ with $\alpha$ verifying \eqref{esta}, $v_\l$ is uniformly bounded, hence $u^*$ is a regular solution if $n$ satisfies \eqref{ndelta}.

\medskip\noindent
{\sl Case 2}: $\d \leq 1$ and $\frac{1}{1 + \d} < \mu \leq \frac{1+\delta}{4\delta}$.

\smallskip
Now we take $\alpha\in \left(\frac 1 2(1-\mu+\epsilon), \frac{1}{2}\right)$, the formulas \eqref{3.5} and \eqref{3.6} imply then
$$\left(1+\frac{\mu-\epsilon}{2\alpha-1}\right)\xi(t)\geq\frac
{\alpha^2}{2\alpha-1}g^{2\alpha-1}(t)g'(t)+C.$$
The inequality \eqref{xi-g} still holds true. Proceeding as for Case 1, we see that for $\d' < \d$ but nearby, $\e > 0$ small and $\beta = 2 - \e$, there exists
$$\alpha \in \left(\frac{1-\mu+\epsilon}2, \frac{\d'}{1+\d'} + \frac{\sqrt{\d'(1+\d')(\mu-\epsilon)-\d'}}{1+\d'}\right) \subset \left(\frac{1-\mu+\epsilon}2, \frac12\right)$$ such that \eqref{da} is satisfied. Hence we conclude exactly as in {\sl Case 1}.\hfill{$\Box$}

\subsection{Proof of Theorem 1.5}
Without loss of generality, assume again $g(0)=1$. Take now
$\phi(t)=te^{\alpha h(t)}$, where $\alpha>0$ is a constant to be
determined, then
\begin{align*}
    \xi(t)& =\int_0^t \left[1+s\alpha h'(s)\right]^2e^{2\alpha h(s)}ds\\
           & =\int_0^t \left[1+ 2s\alpha h'(s)\right] e^{2\alpha h(s)}ds+\int_0^t
           \alpha^2s^2h'^2(s)e^{2\alpha h(s)}ds\\
           & =te^{2\alpha h(t)}+K(t).
\end{align*}
Thus, for $t \geq t_0$,
\begin{align*}
    \frac {2K(t)}{\alpha} =2\alpha\int_0^t s^2h'^2(s)e^{2\alpha
    h(s)}ds &=C+\int_{t_0}^t s^2h'(s)d\left(e^{2\alpha
                         h(s)}\right)\\
                         &\leq C+t^2h'(t)e^{2\alpha
                         h(t)}-\int_{t_0}^t e^{2\alpha
                         h(s)}d\left(s^2h'(s)\right),
\end{align*}
where the last integration is considered in the sense of Stieltjes.
The monotonicity of $s^2h'$ in $[t_0,\infty)$ yields
$$K(t)\leq \frac \alpha 2 t^2h'(t)e^{2\alpha h(t)}+C, \quad \forall\; t \geq t_0.$$
So we get
$$\xi(t)\leq C+ \left[t+\frac \alpha 2 t^2h'(t)\right] e^{2\alpha
h(t)},\quad \forall\; t \geq 0.$$ Using (\ref{psi and xi}) (we
drop the index $\lambda$),
\begin{align*}
    & \int_\Omega e^\gamma\left[e^{h(v)}+vh'(v)e^{h(v)}-ve^{h(v)}-1\right] v^2e^{2\alpha
    h(v)}\\
    \leq & \; \frac2\beta\int_\Omega e^\gamma\left( 1+ve^{h(v)}\right) \xi(v)
    + C \int_\Omega e^\gamma v^2e^{2\alpha h(v)}\\
    \leq & \; \frac2\beta\int_\Omega e^\gamma \left(1+ve^{h(v)}\right)\left[C+ve^{2\alpha h(v)}+\frac \alpha 2 v^2h'(v)e^{2\alpha
    h(v)}\right] + C\int_\Omega e^\gamma v^2e^{2\alpha h(v)},
\end{align*}
By Young's inequality,
\begin{align}
\label{estalpha}
\begin{split}
   & \left(1-\frac\alpha\beta\right)\int_\Omega e^\gamma v^3h'(v)e^{(2\alpha+1)h(v)}\\ \leq & \; C\int_\Omega e^\gamma \left[1+ v^2h'(v)e^{2\alpha h(v)}+v^3e^{(2\alpha+1)h(v)}\right].
\end{split}
\end{align}
Moreover, $g = o(g')$ at infinity yields $\lim_{t \to \infty} h'(t) = \infty$, hence
\begin{align*}
    \frac {t^2h'(t)e^{2\alpha
    h(t)}+t^3e^{(2\alpha+1)h(t)}}{t^3h'(t)e^{(2\alpha+1)h(t)}} = \frac 1 {g(t)-1}+\frac 1 {h'(t)} \rightarrow 0 \; \mbox{ as } \;t \rightarrow \infty.
\end{align*}
Fix $\beta \in (\alpha, 2)$, the inequality \eqref{estalpha} implies
$$\int_\Omega \frac {[g(v)-1]^{2\alpha+1}}{v^{2\alpha}}=\int_\Omega v
e^{(2\alpha+1)h(v)} \leq C + \int_\Omega v^3h'(v)e^{(2\alpha+1)h(v)}\leq
C.$$
Recall that $g$ is superlinear, we obtain $\|g(v)\|_1 \leq C$. Consider again $w$ satisfying $L(w) = \l g(v)$ in $\O$ and $w = 0$ on $\p\O$, as $v \leq w$ in
$\Omega$ by maximum principle,
$$\int_\Omega\frac{(g(v)-1)^{2\alpha+1}}{w^{2\alpha}}\leq
C.$$
Following the proof of Lemma 2.1 in \cite{yz1} (we just need a minor adjustment,
say define $\Omega_1 =\{x\in\Omega: g(v)>w^T\}$ instead, here $T > 0$ is a suitable constant), we can obtain that if $2\alpha+1> n/2$, $w$ is uniformly bounded in $L^\infty(\Omega)$, so does $v$. Taking $2 > \beta > \alpha > 7/4$, the result holds for $n \leq 9$.\hfill{$\Box$}

\subsection{Proof of Theorem 1.6}
Here we choose $\phi(u) = (1 - u)^{-\alpha} - 1$ in \eqref{est1}. For $2\l >
\l^*$ and $\e > 0$,
$$\left(p - \frac{2\alpha^2}{\beta(2\alpha + 1)} - 2\epsilon\right)\int_\O \frac{e^\gamma}{(1 - u)^{p+2\alpha +1}} \leq C, \quad \forall\;\beta \in [1, 2).$$
We have used $f'(u)(1 - u) \geq (p - \epsilon)f(u) - C$ in
$[0, 1)$ by \eqref{H2'}. As $\e > 0$ is arbitrary,
$$\int_\O \frac{1}{(1 - u)^{p+2\alpha +1}} \leq C$$
provided that
$$p > \frac{\alpha^2}{2\alpha + 1},\quad \mbox{i.e. when }\;\alpha < p + \sqrt{p(p+1)}.$$ Therefore $\|(1 -
u)^{-1}\|_q \leq C$ if $q < 1 + 3p + 2\sqrt{p(p+1)}$. For any $\e > 0$, as $f'(u)(1-u) \leq (p + \e)f(u) + C_\e$ in
$[0, 1)$ by \eqref{H2'}, we have $f(u) \leq C(1 - u)^{-p - \e}$, consequently
$$g(v) = e^vf(1 - e^{-v}) = \frac{f(u)}{1 - u} \leq C(1 - u)^{-1 -p - \e},$$ hence $\|g(v)\|_r \leq C$ when
$$r < \frac{1 + 3p + 2\sqrt{p(p+1)}}{p+1+\e}.$$
According to Lemma \ref{comp}, the proof is done by taking $\e \to 0^+$.\hfill{$\Box$}

\section{Radial case}
\setcounter{equation}{0}
\setcounter{lemma}{0}
\setcounter{theorem}{0}
As we have mentioned, when $c = -\nabla\gamma$, the equation $(P_\l)$ is rewritten as
\begin{align}
\label{selfadjoint}
-{\rm div}(e^\gamma\nabla u) = \l e^\gamma f(u).
\end{align}
With the variational structure, the stability of minimal solutions $u_\l$ is equivalent to
\begin{align}
\label{stable2}
\int_\O e^\gamma |\nabla \psi|^2 \geq \l \int_\O e^\gamma f'(u_\l)\psi^2, \quad \forall\; \psi\in H_0^1(\O).
\end{align}
Moreover, for any $C^1$ functions $\phi$ and $\xi$ satisfying
$\phi(0) = \xi(0) = 0$ and $\xi'=\phi'^2$, the estimate \eqref{est1} is replaced by
\begin{align*}
\int_\Omega e^\gamma f'(u_\l)\phi^2(u_\l) \leq
\int_\Omega e^\gamma f(u_\l)\xi(u_\l).
\end{align*}
Taking now $\phi(t) = f(t) - f(0)$ and working as for Theorem 1 in \cite{yz2}, we have
\begin{lemma}
\label{gradientH2}
When $c = \nabla\gamma$, the extremal solution $u^* \in H^2\cap H_0^1(\O)$. More precisely,
\begin{align}
\label{H2}
\int_\O f'(u_\l)f(u_\l) \leq C, \quad  \forall\; \l \in (0, \l^*].
\end{align}
\end{lemma}
When $\O = B_1$ is the unit ball, $\gamma(x) = \gamma(r)$ with $r = |x|$, $u_\l$ is radial by uniqueness of the minimal solution and satisfies
\begin{align}
\label{radial}
-u'' - \frac{n-1}{r}u' - \gamma'u' = \l f(u)\quad \mbox{in } (0, 1],
\end{align}
with $u'(0) = 0$ and $u(1) = 0$. Our main result in this section is the regularity of the extremal solution $u^*$ for any $f$ satisfying $(H)$ provided $n = 2$ and the optimal estimate for $u'$ claimed in Theorem \ref{radialu'}.

\medskip
The method we use is similar to \cite{cc, v}, but the uniform boundedness of $\|u_\l\|_{C^1}$ is not enough to claim the regularity of $u^*$, because a singular $u^*$ could be Lipschitz in many cases (see Remark \ref{1.1}). In fact, the estimate \eqref{H2} is crucial for our proof.

\medskip
As in \cite{cc, v}, since $u_\l'(r) \leq 0$ by maximum principle or equation \eqref{radial}, the boundedness of $\|u_\l\|_{H_0^1}$ implies that for any $k\in \N$, $r > 0$, $\|u_\l\|_{C^k\left(\overline B_1\setminus B_r\right)} \leq C_{k, r}$, $\forall\; \l \in (0, \l^*]$. So we concentrate our attention near the origin. Derivating the equation \eqref{radial} or \eqref{selfadjoint} with respect to $r$,
\begin{align*}
-{\rm div}\left(e^\gamma\nabla u'\right) = e^\gamma u'\left[\l f'(u)- \frac{n-1}{r^2} + \gamma''\right] \; \mbox{ in } (0, 1].
\end{align*}
Using $\psi = r\eta(r) u_\l'(r)$ as test function in \eqref{stable2} with $\eta \in H_0^1(B_1)\cap C(\overline B_1)$, by similar calculation as for Lemma 2.1 in \cite{cc}, we obtain
\begin{align}
\label{radialest1}
\int_{B_1} e^\gamma \Big[|\nabla(r\eta)|^2 - (n-1)\eta^2 + \gamma''r^2\eta^2\Big]u_\l'^2 \geq 0, \quad \forall\; \l \in (0, \l^*].
\end{align}

\subsection{Proof of Theorem 1.1}
For simplicity, we drop the index $\l$. All estimates below hold uniformly for $\l$. First as $u_\l$ is radial, by maximum principle, we see that $u$ is decreasing in $r$. Since $f$ and $f'$ are nondecreasing functions according to $(H)$, the estimate \eqref{H2} implies (as $n = 2$)
\begin{align*}
\pi r^2 f'(u(r))f(u(r)) \leq \int_{B_r} f'(u)f(u) \leq C, \quad \forall\; r \in (0, 1].
\end{align*}
By Lemma \ref{ff'}, we have
\begin{align}
\label{f}
f(u(r)) \leq \frac{C}{r} \quad \mbox{for all } r \in (0, 1].
\end{align}
Let $r_0 \in (0, \frac{1}{2}]$. Let $\eta$ be a radial function in $H_0^1(B_1)\cap C^0(\overline B_1)$ such that
$$\eta(r) = \left\{\begin{array}{ll}
r_0^{-1} & \mbox{if } r < r_0;\\
r^{-1} & \mbox{if } r_0 \leq r \leq \frac{1}{2},
\end{array}
\right.$$
and $\eta$ be a fixed $C^1$ function in $\overline B_1\setminus B_{1/2}$, independent of $r_0$. The direct calculation yields
$$|\nabla(r\eta)|^2 - \eta^2 + \gamma''r^2\eta^2 = \left\{\begin{array}{ll}
\gamma''r^2r_0^{-2} & \mbox{if } r < r_0;\\
\gamma'' - r^{-2} & \mbox{if } r_0 < r \leq \frac{1}{2}.
\end{array}
\right.$$
Using \eqref{radialest1}, as $u$ is uniformly bounded in $H^1(B_1)$ by Proposition \ref{H01} and $r^2r_0^{-2} \leq 1$ in $[0, r_0]$, we get
$$\int_{r_0}^{\frac{1}{2}} \frac{u'(r)^2}{r} dr \leq C.$$
Tending $r_0$ to $0$, there holds
\begin{align}
\label{radialest2}
\int_0^1 \frac{u'(r)^2}{r} dr \leq C.
\end{align}
Consider the following test function used in \cite{v}: For any $r \leq \frac{1}{2}$ and $0 < r_0 < r$,
$$\eta(s) = \left\{\begin{array}{ll}
(rr_0)^{-1} & \mbox{if } s < r_0;\\
(rs)^{-1} & \mbox{if } r_0 \leq s < r;\\
s^{-2} & \mbox{if } r \leq s \leq \frac{1}{2}.
\end{array}
\right.$$
Applying again \eqref{radialest1} and combining with \eqref{radialest2}, we obtain finally (with $r_0 \rightarrow 0$)
\begin{align}
\label{radialest3}
\int_0^r \frac{u'(s)^2}{s} ds \leq Cr^2, \qquad \forall\; r \leq 1.
\end{align}
As $\left(e^\gamma ru'\right)' = -\l e^\gamma rf(u)$ with $n = 2$, so $e^\gamma ru'$ is nonincreasing in $r$. Then $u'(s) \leq Cru'(r)/s$ for $s \in [r, 1]$, hence $u'(s) \leq Cu'(r) \leq 0$ for any $s\in [r, 2r]$ if $r \leq \frac{1}{2}$. By \eqref{radialest3}, for any $0 < r \leq \frac{1}{2}$,
\begin{align*}
C_1r^2 \geq \int_0^{2r} \frac{u'(s)^2}{s} ds \geq \int_r^{2r} \frac{u'(s)^2}{s} ds \geq \frac{C_2}{r}\int_r^{2r} u'(r)^2ds = C_3u'(r)^2.
\end{align*}
That means
\begin{align}
\label{estu'}
|u'(r)| \leq Cr \quad \mbox{in } \;[0, 1].
\end{align}
However, we need to consider also $u''(r)$ as explained above. Let
$$G(r) = e^\gamma r u' \quad \mbox{and}\quad \Psi(r) = -2G(\sqrt{r}) - M\int_0^r (r-s)f\left(u(\sqrt{s})\right)ds$$
where $M$ is a constant to be chosen. Using $G' = - \l e^\gamma rf(u)$,
\begin{align*}
\Psi''(r) & = \left[\l e^{\gamma(s)}f'\left(u(s)\right)\frac{u'(s)}{2s} + \l e^{\gamma(s)} f\left(u(s)\right)\frac{\gamma'(s)}{2s} - Mf\left(u(s)\right)\right]\Big|_{s = \sqrt{r}}\\
& \leq \left[\l e^{\gamma(s)} f\left(u(s)\right)\frac{\gamma'(s)}{2s} - Mf\left(u(s)\right)\right]\Big|_{s = \sqrt{r}}\\
& \leq C_0f\left(u(\sqrt{r})\right) - M f\left(u(\sqrt{r})\right).
\end{align*}
For the last line, we used $|\gamma'(s)|/s \leq C$ in $[0, 1]$ since $\gamma$ is a smooth function (so $\gamma'(0) = 0$). Fix $M > C_0 + 1$, $\Psi$ is then concave in $[0, 1]$. On the other hand, by \eqref{f}
\begin{align*}
\Psi'(r) = \l e^{\gamma(\sqrt{r})}f\left(u(\sqrt{r})\right) - M\int_0^r f\left(u(\sqrt{s})\right) ds \geq C\l f(0) - CM\sqrt{r}.
\end{align*}
There exists $r_1 > 0$ small enough such that $\Psi' \geq 0$ in $[0, r_1]$ with $\l \geq \frac{\l^*}{2}$. Using \eqref{radial}, \eqref{f} and \eqref{estu'}, for $\l \geq \frac{\l^*}{2}$ and $r \leq r_1$,
\begin{align*}
& -e^{\gamma(\sqrt{r})}\left[u''(\sqrt{r}) + \frac{u'(\sqrt{r})}{\sqrt{r}} + \gamma'u'(\sqrt{r})\right] - CM\sqrt{r}\\ \leq & \; \Psi'(r) \leq \frac{\Psi(r)}{r} \leq -2e^{\gamma(\sqrt{r})}\frac{u'(\sqrt{r})}{\sqrt{r}} \leq C.
\end{align*}
Applying one more time \eqref{estu'}, we see that $u''(\sqrt{r}) \geq -C$ for any $\l \geq \frac{\l^*}{2}$ and $r \leq r_1$. Otherwise, by \eqref{radial} and \eqref{estu'}, $u''(r) \leq -u'(r)r^{-1} - \gamma'(r)u'(r) \leq C$, we claim then
\begin{align*}
\|u''\|_\infty \leq C, \quad \forall\; \l \geq \frac{\l^*}{2}.
\end{align*}
Combining with \eqref{radial} and \eqref{estu'}, it means $\|\l f(u)\|_\infty \leq C$, no singularity will occur.\hfill{$\Box$}

\subsection{Proof of Theorem 1.2}
As above, we drop the index $\l$ and all estimations hold uniformly for $\l$. First, repeating the proof of Theorem 1.8, c) in \cite{cc}, we obtain $f'(u(r)) \leq Cr^{-2}$ in $(0, 1]$. Using Lemma \ref{ff'} with \eqref{radialest1}, $f(u(r)) \leq Cr^{-2}$ in $(0, 1]$. Consequently, by \eqref{radial}, for $n \geq 3$,
\begin{align*}
0 \leq -e^\gamma r^{n-1}u'(r) = \int_0^r e^{\gamma(s)} s^{n-1}f(u(s))ds \leq C\int_0^r s^{n-3}ds \leq Cr^{n-2}.
\end{align*}
Hence
\begin{align}
\label{generalu'}
|u'(r)| \leq \frac{C}{r}.
\end{align}
Let $\eta$ be a radial function in $H_0^1(B_1)\cap C^0(\overline B_1)$ such that
$$\eta(r) = \left\{\begin{array}{ll}
r_0^{-\sqrt{n-1}} & \mbox{if } r < r_0;\\
r^{-\sqrt{n-1}} & \mbox{if } r_0 \leq r \leq r_1.
\end{array}
\right.$$
in $\overline B_{r_1}$ and be a fixed $C^1$ function in $\overline B_1\setminus B_{r_1}$, here $r_0$ is any constant in $(0, r_1)$, $r_1 > 0$ is a small constant to be determined. Therefore
$$
|\nabla(r\eta)|^2 - (n-1)\eta^2 + \gamma''r^2\eta^2 = \left\{\begin{array}{ll}
\left(\gamma''r^2 + 2 - n\right)r_0^{-2\sqrt{n-1}} & \mbox{if } r < r_0;\\
\left(\gamma''r^2 - 2\sqrt{n-1} + 1\right) r^{-2\sqrt{n-1}} & \mbox{if }r \in [r_0, r_1].
\end{array}
\right.
$$
We fix $r_1 > 0$ small enough such that
$$\max_{r \in [0, r_1]}\left\{\gamma''r^2 \right\} < \min\left(n - 2, 2\sqrt{n-1} - 1\right).$$
By \eqref{radialest1}, as $|\nabla(r\eta)|^2 - (n-1)\eta^2 + \gamma''r^2\eta^2 \leq 0$ for $r\in [0, r_0]$,
\begin{align*}
\int_{r_0}^{r_1} u'^2(r) r^{n - 1-2\sqrt{n-1}}dr \leq C.
\end{align*}
Tending $r_0$ to $0$, we have
\begin{align}
\label{radialest4}
\int_0^{r_1} u'^2(r) r^{n - 1-2\sqrt{n-1}}dr \leq C.
\end{align}
Now we take another test function used in \cite{v},
$$\eta(r) = \left\{\begin{array}{ll}
r_0^{-\sqrt{n-1}-1} & \mbox{if } r < r_0;\\
r^{-\sqrt{n-1}-1} & \mbox{if } r_0 \leq r \leq r_1.
\end{array}
\right.$$
Combining \eqref{radialest1} and \eqref{radialest4}, we conclude then
\begin{align*}
\int_0^{r_0} u'^2(r) r^{n - 1}dr \leq Cr_0^{2+ 2\sqrt{n-1}}, \quad \forall\; r_0 \in [0, r_1].
\end{align*}
By the monotonicity of $e^\gamma r^{n-1}u'$, similarly as for \eqref{estu'}, it holds
\begin{align*}
|u'(r)| \leq C r^{-\frac{n}{2}+1 + \sqrt{n-1}}, \quad \forall\; r \in [0, 1].
\end{align*}
Finally, combining with \eqref{generalu'}, we are done (in fact, $-\frac{n}{2}+1 + \sqrt{n-1} \leq -1$ for $n \geq 10$). \hfill{$\Box$}

\bigskip
\noindent {\bf Acknowledgments} {\sl Part of the work was completed
during X.L.'s visit to University of Connecticut (Uconn) with the
financial support of CSC. She would like to thank the Department of
Mathematics of Uconn for its warm hospitality. She thanks also Prof. Ryzhik for useful discussion. D.Y. is supported by
the French ANR project referenced ANR-08-BLAN-0335-01. F.Z. is
supported in part by NSFC No. 10971067, the ``basic research project
of China" No. 2006CB805902 and Shanghai project 09XD1401600.}

\end{document}